\documentclass[12pt,leqno]{article}

\usepackage{amssymb, amsfonts, eucal,amsmath}

\def\li {\limits }
\def\bbr {\Bigl)}
\def\bbl {\Bigl(}
\def\br {\bigl)}
\def\bl {\bigl( }
\def\om {\omega }

\def\Ga {\Gamma}

\def\ol {\overline }
\def\ul {\underline }

\def\oh {\ol{H}}
\def\ooh {\ol{\ol{H}}}
\def\lo {\left }

\def\rr {\right }

\def\wti {\widetilde }
\def\wh {\widehat }

\def\jn {j(\nu) }

\def\ti {\tilde }

\newtheorem{theorem}{Theorem}
\newtheorem{lemma}[theorem]{Lemma}
\newtheorem{remark}{Remark}

\def\proof{\medskip\noindent{\bf Proof.} }
\def\eop{\hfill$\square$ \\ \mbox{} }

\newcommand{\sect}[1]{\section{#1} \setcounter{equation}{0} }

\begin{document}

\begin{center}

\renewcommand{\thefootnote}{\empty }

{\LARGE Jackson's type estimate  \\ of nearly coconvex approximation}



{\large G.A.~Dzyubenko}

\vspace{0.2cm}

{\small\it Yu.A.Mitropolsky International Mathematical Center of NASU, \\
Tereschenkivska str. 3, 01601 Kyiv-4, Ukraine \\
dzyuben@gmail.com}

\end{center}

\vspace{0.2cm}

\begin{abstract} Suppose that a continuous on the real axis $2\pi$-periodic function $f$
changes its convexity at $2s,\ s\in\Bbb N,$ points $y_i$ on each period:
$-\pi\le y_{2s}<y_{2s-1}<...<y_1<\pi,$ and for the rest $i\in\Bbb Z,$ the points $y_i$ are defined periodically. In the paper, for each $n\ge N,$ a trigonometric polynomial $P_n$ of order $cn$ is found such that: $P_n$ has the same convexity as $f,$ everywhere except, perhaps, the small neighborhoods of the $y_i:$
$$
(y_i-\pi/n,y_i+\pi/n)
$$
and
$$
\|f-P_n\|\le c(s)\,\omega_4(f,\pi/n),
$$
where $N$ is a constant depending only on $\min\limits_{i=1,...,2s}\{y_i-y_{i+1}\},\ c$ and $c(s)$ are  constants depending only on $s,\ \omega_4(f,\cdot)$ is the modulus of continuity of the $4$-th order of the function $f,$ and $\|\cdot\|$ is the max-norm.
\end{abstract}

\vspace{0.2cm}

\sect{Introduction and the main theorem}

By $C$ we denote the space of continuous $2\pi$-periodic functions
$f:\Bbb R\rightarrow\Bbb R$ with the uniform norm
$$
\left \Vert  f\right \Vert =\max\limits_{x\in\Bbb R}\left
|f(x)\right |,
$$
and by $\Bbb T_n,\ n\in\Bbb N,$ denote the space of trigonometric polynomials
$$
P_n(x)=a_0+\sum_{j=1}^{n}(a_j\cos
jx+b_j\sin jx),\quad  a_j\in\Bbb R, \
b_j\in\Bbb R,
$$
of degree $ \le n.$ Recall the classical Jackson-Zygmund-Akhiezer-Stechkin estimate
(obtained by Jackson for $k=1,$ Zygmund and Akhiezer for $k=2,$ and Stechkin for $k\ge
3,\ k\in\Bbb N$): {\it if a function $f\in C,$ then for each $n\in\Bbb N$ there is a polynomial $P_n\in\Bbb T_n$ such that
\begin{equation}\label{1}
\|f-P_n\|\le c(k)\, \omega_k \left(f,\pi/n\right),
\end{equation}
where $c(k)$ is a constant depending only on $k,$ and
$\omega_k \left(f,\cdot\right)$ is the modulus of continuity of order $k$ of the function $f.$} For details, see, for example,
\cite
{dzay_shev}.

In 1968 Lorentz and Zeller \cite{LZ} for $k=1$ obtained a bell-shaped analogue of the inequality
\eqref{1}, {\em i.e., } when bell-shaped (even and nonincreasing on $[0,\pi]$) $2\pi$-periodic functions are approximated by bell-shaped polynomials.

In papers \cite{petay} and \cite{Zalizko} two coconvex analogues of the inequality \eqref{1} were proved for $k=2$ and $k=3,$ respectively. Moreover, in \cite{ZalizkoContr} arguments from the papers \cite{shved_otr1}, \cite{shved_otr2} of Shvedov and \cite{devlevshev} of DeVore, Leviatan and Shevchuk were used to show that for $k>3$ there is no coconvex analogue of the inequality \eqref{1}.

Nevertheless, as we know from the coconvex approximation on a closed interval (by algebraic polynomials, see, for details \cite{DzGi2006}) {\it if some relaxation of the condition of coconvexity for the approximating polynomials is allowed, then an extra order of the approximation can be achieved}, and, as it seems, no more than one extra order, though the corresponding counterexample is not constructed yet.

So, in the paper in Theorem 1 we prove a trigonometric analogue of the algebraic result \cite{DzGi2006}.
 To write it we give necessary notations.

 Suppose that on $[-\pi ,\pi)$ there are $2s,\ s\in\Bbb N,$ fixed points
$y_i:$
$$
-\pi \le y_{2s}<y_{2s-1}<\dots<y_1<\pi ,
$$
while for other indices $i\in\Bbb Z,$ the points $y_i$ are defined periodically by the equality
$$
y_i=y_{i+2s}+2\pi \quad ({\rm i.e.,} \ y_0=y_{2s}+2\pi ,...,
y_{2s+1}=y_1-2\pi ,... ).
$$
Denote $ Y:=\{ y_i\}_{i\in\Bbb Z}. $
By $\Delta^{(2)}(Y)$  we denote the set of all functions $f\in C$ which are convex on
$[y_1,y_0],$  concave on $[y_2,y_1],$
convex on $[y_3,y_2],$ and so on. The functions in $\Delta^{(2)}(Y_s)$ are
{\it coconvex} with one another.
Note, if a function $f$ is twice differentiable, then
$f\in\Delta^{(2)}(Y)$ if and only if
$$
f''(x)\Pi(x)\ge 0,\ x\in\Bbb R,
$$
where
$$
\Pi(x):=\Pi(x,Y):=\prod_{i=1}^{2s}\sin\frac{x-y_i}{2}\quad\quad
(\Pi(x)>0,\ x\in(y_1,y_0)).
$$
\begin{theorem}\label{thm1}
 If a function $f\in\Delta^{(2)}(Y),$ then there exists a constant $N(Y)$ depending only on $
\min\limits_{i=1,...,2s}\,\{y_{i}-y_{i+1}\}$ such that for each $n\ge N(Y)$ there is a polynomial
$P_n\in\Bbb T_{cn}$ satisfying
\begin{equation}\label{nearly}
P_n''(x)\Pi(x)\ge 0,\quad x\in\Bbb R\setminus \cup_{i\in\Bbb Z}\left(y_i-\pi/n,y_i+\pi/n\right),
\end{equation}
\begin{equation}\label{estimate}
\|f-P_n\|\le c(s)\,\omega_4(f,\pi/n),
\end{equation}
where $c$ and $c(s)$ are constants depending only on $s.$
\end{theorem}
The following Theorem \ref{thm2} is a simple corollary of Theorem 1 and Whitney's inequality
 \cite{whit}
$$
\|f-f(0)\|\le 3\,\omega_4(f,4\pi).
$$
\begin{theorem}\label{thm2}
 If a function $f\in\Delta^{(2)}(Y),$ then for each $n\in\Bbb N$ there is a polynomial
$P_n\in\Bbb T_n$ such that
\begin{equation}\label{nearly2}
P_n'(x)\Pi(x)\ge 0,\quad x\in\Bbb R\setminus \cup_{i\in\Bbb Z}\left(y_i-c/n,y_i+c/n\right),
\end{equation}
\begin{equation}\label{estimate2}
\|f-P_n\|\le C(Y)\,\omega_4(f,\pi/n),
\end{equation}
where $c$ is a constant depending only on $s,$ and $C(Y)$ is a constant depending only on $
\min\limits_{i=1,...,2s}\,\{y_{i}-y_{i+1}\}.$
\end{theorem}

\begin{remark}
We believe that $\omega_4$ in \eqref{estimate} and \eqref{estimate2}  cannot be replaced by
$\omega_k$ with $k>4.$ Also we believe that the constants $N(Y)$ and $C(Y)$ in Theorems 1 and 2 cannot be replaced by constants independent of $ \min\limits_{i=1,...,2s}\,\{y_{i}-y_{i+1}\}$ (and depending, say, on $s$). These both assumptions are not made further in the paper. Also, we do not pay attention to the constant $c$ in the both theorems, i.e., we did not try to replace it by an absolute  constant or/and by a smallest possible one.
\end{remark}

\sect{Auxiliary facts I}

For each $n\in\Bbb N$ denote
$$
h:=h_n:=\frac{\pi}{n},\quad x_j:=x_{j,n}:=-j\,h,\quad
I_j:=I_{j,n}:=[x_j,x_{j-1}],\quad
j\in\Bbb Z.
$$
Let
$
m\in\{1,2,3,10,20,30\}.
$
For a fixed  $Y=\{y_i\}_{i\in\Bbb Z}$ and a fixed $n$
denote
$$
O_{i,m}:=O_{i}(Y,n,m):=(x_{j+m+1},x_{j-m})\quad\text{if}\quad
y_i\in[x_j,x_{j-1})=:[x_{j_i},x_{j_i-1}).
$$
Set
$$
O_m:=O(Y,n,m):=\mathop{\bigcup}\limits_{i\in\Bbb Z}O_{i,m}.
$$
We will write
$$j\in H(Y,n,m)\quad \text{if}\quad x_j\subset\Bbb R\setminus O_m.$$
Let
$$
H_m:=\left\{j:\,j\in H(Y,n,m),\ |j|\le n\right\}.
$$

Choose $N(Y):=N(Y,30)\in\Bbb N$ sufficiently large so that
\begin{equation}\label{OO}
 O_{i,30}\cap O_{i-1,30}=\emptyset
\end{equation}
for all $n\ge N(Y)$ and all $i=1,...,2s$ (thus, $N(Y)$
 depends only on  $\min\limits_{i=1,...,2s}\,\{y_{i}-y_{i+1}\}$). In what follows $n>N(Y).$

Denote
$$
\chi(x,a):=\left\{
\begin{array}{lll}
0,  & \mbox{if} &  x\le a, \\
1,  & \mbox{if} & x>a,
\end{array}
\right. \ a\in\Bbb R,\quad \chi_j(x):=\chi(x,x_j),\quad (x-x_j)_+:=(x-x_j)\chi_j(x),
$$
$$
\Gamma_j(x):=\Gamma_{j,n}(x):=\min\left\{1,\frac{1}{n\,\left|\sin
\frac{x-(x_j+h/2)}{2}\right|}\right\},\quad j\in\Bbb Z,\ n\in\Bbb
N,
$$
and note that
\begin{equation}\label{2.2}
\left\Vert\sum\limits_{j=1-n}^{n}\Gamma_j^2\right\Vert < 6,
\end{equation}
for details, see \cite{plesh}.

For each $j\in\Bbb Z$ and $b\in \Bbb N$ we set the positive polynomial $J_j\in\Bbb T_{(n-1)b},\ n\in\Bbb N,$
\begin{equation}\label{Jj}
J_j(x):=J_{j,n}(x):=\left(\frac{\sin\frac{n(x-x_j)}{2}}{\sin\frac{x-x_j}{2}}\right)^{2b}+
\left(\frac{\sin\frac{n(x-x_{j-1})}{2}}{\sin\frac{x-x_{j-1}}{2}}\right)^{2b}
\end{equation}
(i.e., the sum of two "adjacent"  kernels of Jackson type).

For each $j\in H_{10}$ denote
\begin{equation}\label{tj}
t_j(x):=t_{j,n}(x,b,Y):=\frac{\int_{x_j-\pi}^x
J_j(u)\Pi(u)du}{\int_{x_j-\pi}^{x_j+\pi } J_j(u)\Pi(u)du}.
\end{equation}

In what follows $ c_i=c_i(b)=c_i(s,b),\ i=1,...,8,$ stand for positive constants which may depend only on $s$ and $b.$

\begin{lemma}{\rm \cite[Lemma 1]{dp}}. \label{lem1} If $j\in  H_{10}$ and $b\ge
s+2,$ then
\begin{equation}\label{tj'pipi}
t_j'(x)\,\Pi(x)\,\Pi(x_j)\ge 0,\quad x\in\Bbb R,
\end{equation}
\begin{equation}\label{chi-tj}
\left|\chi_j(x)-t_j(x)\right|\le c_1\left(\Gamma_j(x)\right)^{2b-2s-1},
\quad x\in [x_j-\pi ,x_j+\pi],
\end{equation}
\begin{equation}\label{tj'<}
\left|t_j'(x)\right|\le
c_2\frac{1}{h}\,\left(\Gamma_j(x)\right)^{2b-2s},\quad x\in\Bbb R,
\end{equation}
\begin{equation}\label{tj'>}
\left|t_j'(x)\right|\ge
c_3\frac{1}{h}\,\left(\Gamma_j(x)\right)^{2b+2s},\quad x\in\Bbb
R \setminus O_{10},
\end{equation}
\begin{equation}\label{tj'>o}
\left|t_j'(x)\right|\ge
c_3\frac{1}{h}\,\left(\Gamma_j(x)\right)^{2b+2s}\left|\frac{x-y_i}{x_j-y_i}\right|,\quad
x\in O_{i,10},\ i\in\Bbb Z.
\end{equation}
\end{lemma}

Note that Lemma \ref{lem1} is proved by using the inequalities
\begin{equation}\label{2.2''}
\begin{array}{rcl}
\frac{1}{c_4h}\Gamma_j^{2b}(x)\left|\frac{\Pi(x)}{\Pi(x_j)}\right|
&\le &\left|t_j'(x)\right|\ \le \  \frac{c_4}{h}\Gamma_j^{2b}(x)\left|\frac{\Pi(x)}{\Pi(x_j)}\right|,\\
\left|\frac{\Pi(x)}{\Pi(x_j)}\right|\
&\le &2^{2s}\Gamma_j^{-2s}(x), \quad
x\in\Bbb R,\quad j\in H_m,\ m\ge 10,
\end{array}
\end{equation}
\begin{equation}\label{2.2'}
\begin{array}{rcl}
\left|\int_{x}^{x_j+\pi} \Gamma_j^{b}(u) du\right|&\le & c_5\,h\,\Gamma_j^{b-1}(x),
\quad b\in\Bbb N,\quad
x\in[x_j,x_j+2\pi],\\
\left|\int_{x}^{x_j-\pi} \Gamma_j^{b}(u) du\right|\ &\le & c_5\,h\,\Gamma_j^{b-1}(x),
\quad b\in\Bbb N,\quad
x\in[x_j-2\pi,x_j],
\end{array}
\end{equation}
for details, see \cite{plesh}.

For each $j\in H_{20}$ set the function
\begin{equation}\label{tauj}
\tau_j(x):=\tau_{j,n}(x,b,t_j):=\alpha\int_{x_j-\pi}^{x}t_{j+10}(u)du
+(1-\alpha)\int_{x_j-\pi}^{x}t_{j-10}(u)du,
\end{equation}
where the number $\alpha\in [0,1]$ is chosen from the condition
$$
\tau_{j}(x_j+\pi)=\pi
$$
(note that the inequalities $0\le\alpha\le 1$ follow from the estimate \eqref{chi-tj} and the choice of the indices $j\pm 10$ if $b\ge s+2,$ for details, see \cite[p. 923]{petay}).

Note that the functions  $t_j$ and $\tau_j$ can be expressed on $\Bbb
R$ as
\begin{equation}\label{tjR}
t_j(x)=\frac{1}{2\pi }x+\hat R_j(x),\quad j\in  H_{10},
\end{equation}
\begin{equation}\label{taujR}
\tau_j(x)=\frac{1}{4\pi }x^2+\frac{\pi -x_j}{2\pi}x+\tilde R_j(x),
\quad j\in H_{20},
\end{equation}
where $\hat R_j$ and $\tilde R_j$ are polynomials from $\Bbb T_{ c_6n}
$ (see similar cases in \cite{plesh} and \cite{petay}, respectively).

In what follows $c>0$ denote different absolute constants or constants depending only on
 $s.$ They can be different even if they are in the same line.

 Denote two functions $\wti{t}_j$ and ${\wti{\tau}}_j:$
$$
{\wti{t}}_j(x):={\wti{t}}_{j,n}(x,b):=
\bar t_{j}(x)+\sum\limits_{i=1}^{2s}\frac{\chi_j(y_i)-
\bar t_{j}(y_i)}{\hat t_{j_i}(y_i)}\hat t_{j_i}(x),\quad j\in H_{10},
$$
where $\bar t_j(x):=t_{j,n}(x,\bar b,\emptyset) $ is the function defined by \eqref{tj} with
$\Pi(x):\equiv 1$ and $\bar b=b+3,$ and
$$
\hat t_{j_i}(x):=\left(\bar t_{j_i+10}(x)-\breve t_{j_i-10}(x)\right)\frac{\Pi(x,Y_i)}{\Pi(x_{j_i},Y_i)}
$$
is the polynomial, where $j_i$ is an index $j$ such that
$y_i\in [x_j,x_{j-1}),\ i=1,...,2s,\ \breve t_j(x):=t_{j,n}(x,\bar b, \breve Y_i)$
is the function \eqref{tj} with
$
\breve Y_i:=\left\{y_i-\pi\nu\right\}_{\nu\in\Bbb Z},
$
and
$$
Y_i:=\left(Y\setminus\{y_i+2\pi\nu\}_{\nu\in\Bbb
Z}\right)\cup\{y_{i}^*+2\pi\nu\}_{\nu\in\Bbb Z},
$$
where $ y_i^*$ is the left endpoint of the interval $O_{i,20},$ if $i$ is odd, and -- the right one, if $i$ is even; and
$$
{\wti{\tau}}_j(x):=
{\wti{\tau}}_{j,n}(x,b):=\tau_{j,n}(x,b,\bar t_j)
+\sum\limits_{i=1}^{2s}\frac{(y_i-x_j)_+-
\tau_{j,n}(y_i,b,\bar t_j)}{\hat t_{j_i}(y_i)}\hat t_{j_i}(x), \quad j\in H_{20}.
$$
Note that the following Lemma \ref{lem2} can be proved with the arguments similar to \cite[Lemma 5.3]{DzGiSh}.

\begin{lemma} {\rm \cite[Lemmas 4 and 5]{dz_comono}}. \label{lem2}
For each $j\in H_{10}$ and $b\ge 3s+2$ the function
${\wti{t}}_j$
satisfies the relations \eqref{chi-tj}, \eqref{tjR}, and in addition,
\begin{equation}\label{chi-t}
\left|\chi_j(x)-{\wti{t}}_j(x)\right|\le c_7\left(\Gamma_j(x)\right)^{2b-2s-1}\left|\frac{x-y_i}{x_j-y_i}\right|,\quad
x\in O_{i,10},\ i=1,...,2s,
\end{equation}
(in particular, $\chi_j(y_i)-{\wti{t}}_j(y_i)=0$).
For each $j\in H_{20}$ and $b\ge 3s+2$ the function ${\wti{\tau}}_j$
satisfies the relation \eqref{taujR}, and in addition,
\begin{equation}\label{x-tau}
\left|(x-x_j)_+-{\wti{\tau}}_j(x)\right|\le
c_8h\left(\Gamma_j(x)\right)^{2(b-s-1)},\quad x\in [x_j-\pi ,x_j+\pi],
\end{equation}
\begin{equation}\label{x-tauo}
\left|(x-x_j)_+-{\wti{\tau}}_j(x)\right|\le
c_8h\left(\Gamma_j(x)\right)^{2(b-s-1)}\left|\frac{x-y_i}{x_j-y_i}\right|,\quad x\in O_{i,10},\ i=1,...,2s,
\end{equation}
(in particular,
$(y_i-x_j)_+-{\wti{\tau}}_j(y_i)=0$).
\end{lemma}

\sect{Auxiliary facts II}

Since we prove Theorem 1 using an intermediate approximation by a spline, i.e., the inequality
$||f-S+S-P_n||\le ||f-S||+||S-P_n||,$ we describe  the  $S$ in this section.

\noindent
Without special references we will use Whitney inequality \cite{whit}
$$
\left\vert f(x)-L_3(x;a;f)\right\vert \le \omega_{4}
\left(f,(b-a)/{4},[a,b]\right),\quad x\in [a,b],
$$
where $L_3$ is Lagrange polynomial interpolating $f$ at $a,\
a+\frac{b-a}{3},\ b-\frac{b-a}{3}$ and $b.$
Fix $j\in\Bbb Z.$ Let
$$
\Psi_3(x,x_{j}):=(x-x_j)_+(x-x_{j-1})(x-x_{j-2}),\quad d_j:=x_{j-1},
$$
$$
a_\nu :=a_{j,\nu}:=x_{j}\lor x_{j-1} \lor x_{j-2},\quad
\wti h_{\nu}:=-h \lor 0 \lor h,\quad
\wh h_{\nu}:=2h^2\lor -h^2 \lor 2h^2,
$$
if $\nu =1\lor 2\lor 3$ respectively. In the following
$\nu\in\{1,2,3\}$ only.

Introduce three functions $\Psi_{j,\nu}\in\Bbb C$ coinciding with
$\Psi_3(x,x_j)$ almost everywhere
$$
\Psi_{j,\nu}(x):=\Psi_3(x,x_j)\,\chi (x, a_\nu
)=(x-a_\nu)^3_++3\,\wti h_{\nu}(x-a_\nu)^2_++\wh
h_{\nu}(x-a_\nu)_+.
$$
That is,
\begin{equation}\label{1213}
\Psi_{j,\nu}(x)=\Psi_3(x,x_j),\ x\in\Bbb R\setminus \lo[x_j,a_\nu
\rr];
\quad
\lo|\Psi_3(x,x_j)-\Psi_{j,\nu}(x)\rr|\le c\,h^3,\ x\in
[x_{j},a_{\nu}],
\end{equation}
\begin{equation}\label{14}
\Psi_{j,\nu}(x)=\int\limits_{d_j-\pi}^x \left( 6\int\limits_{d_j-\pi}^t \left((u-a_\nu
)_++\wti h_{\nu}\chi(u,a_\nu )\right)du +\wh h_{\nu }\chi (t, a_\nu )\right)dt,
\end{equation}
and for $\nu_1,\ \nu_2\in\{1,2,3\},$ we have
\begin{equation}\label{17}
\frac{\Psi_{j,\nu_1}''(x)-\Psi_{j-1,\nu_2}''(x)}{3h}=
\frac{6(x-d_j)-6(x-d_{j-1})}{3h}=2,\quad x\in
(\max\{a_{\nu_1},a_{\nu_2}\},\infty ).
\end{equation}

Without loss of generality suppose that $y_{1}=x_{30}$ (i.e., points $Y$ are far from $-\pi$ and $\pi$), also recall that
$H_3\subset H_2\subset H_1$.

\vskip 0.3cm

\centerline{\it Construction of the nearly coconvex cubic spline}

\medskip

Denote two divided differences of $f$
$$
\aligned
F_j:=[&x_j,x_{j-1},x_{j-2};f], \quad \quad \quad &j&=2-n,...,n,\\
\Phi_j:=[x_{j+1},&x_j,x_{j-1},x_{j-2},x_{j-3};f],\quad  &j&=3-n,...,n-1.
\endaligned
$$
Remark,
$
\Phi_j\,4h=\frac{F_{j+1}-F_{j}}{3h}-
\frac{F_{j}-F_{j-1}}{3h}.
$

Introduce new functions $\Psi_j(x),\ j=3-n,...,n-1.$  For each
$j\in H_2,$ put
$$
\Psi_j(x):=\Psi_{j,2}(x)\quad\text{if}\quad \Phi_j\,\Pi(x_j)\le
0,
\leqno(d.0)
$$
{\it otherwise} put
$$
\hskip-1.4cm\aligned (d.1) \\ (d.2) \\ (d.3) \endaligned
\quad\quad\quad \Psi_j(x):=\lo\{ \aligned
\Psi_{j,1}(x)&\quad\text{if}\quad
|F_{j+1}|>|F_{j}|\ge |F_{j-1}|, \\
\Psi_{j,3}(x)&\quad\text{if}\quad |F_{j+1}|\le |F_{j}|<
|F_{j-1}|,\\
\alpha_j\,\Psi_{j,1}(x)+(1-\alpha_j)\,\Psi_{j,3}(x)&\quad\text{if}\quad
|F_{j+1}|>|F_{j}|< |F_{j-1}|,
\endaligned
\rr .
$$
where
$
\alpha_j:=\frac{F_{j+1}}{F_{j+1}+F_{j-1}}\in (0,1).
$
For other $j=3-n,...,n-1,$ such that $j\notin H_2$ (i.e., $j:\ x_j\in O_{i,2},\ i=1,...,2s,$) put
$$
\Psi_j(x):=\lo\{ \aligned
\Psi_{j,2}(x)&\quad\text{if}\quad\Phi_j\,\Pi(x_j,\ti Y_i\})\le 0,\\
\Psi_{j,1}(x)&\quad\text{otherwise,}
\endaligned \rr .\quad \ti Y_i:=\lo(Y\setminus\{y_i\}\rr)\cup\{x_{j_i+5}\}.
\leqno(d.4)
$$
Set
$$
\Psi_n(x):=\Psi_{3}(x,x_n)
\quad
\Psi_{2-n}(x):
\equiv 0.
\leqno(d.5)
$$

\begin{remark} \label{rem2}
 For the both "strange" cases in (d.4) it is
sufficient to take simply $ \Psi_j(x)=\Psi_{j,2}(x)$  to have the
{\it nearly} coconvexity of the spline below with $f$ however the setting
(d.4) is more convenient to verify the nearly copositivity of
$P_n''$ feather.
\end{remark}

Show that the cubic spline
\begin{equation}\label{26}
S(x)=L_3(x;x_n;f)+4h\sum\limits_{j=3-n}^{n-1}\Phi_j\,
\Psi_j(x),
\end{equation}
or equivalently,
\begin{equation}\label{27}
\aligned
S(x)=L_1(x;x_n;f)+&F_n\lo((x-x_n)(x-x_{n-1})
-\frac{\Psi_n(x)-\Psi_{n-1}(x)}{3h}\rr)\\
+&\sum\limits_{j=3-n}^{n-1}
F_j\,A_j(x)+F_2\frac{\Psi_3(x)}{3h},
\endaligned
\end{equation}
where
$$
A_j(x):=\ol{A}_j(x)-\ul{A}_j(x):=\frac{\Psi_{j+1}(x)-\Psi_j(x)}{3h}-
\frac{\Psi_j(x)-\Psi_{j-1}(x)}{3h},
$$
(having been continued periodically) is nearly coconvex with $f,$ i.e.,
\begin{equation}\label{s''pi}
S''(x)\Pi(x)\ge 0, \ x\in I_j,\ j\in H_3,
 \end{equation}
and satisfies the inequality
\begin{equation}\label{f-s}
\|f-S\|=\|f-S\|_{[-\pi,\pi]}\le c\,\omega_4(f,h)
\end{equation}
(it is convenient to look at the sums in \eqref{26} and \eqref{27} starting from the last addend, for the details
of such a kind of representations, see \cite[Proposition 1]{DzGi}).

With the help of \eqref{26} and \eqref{27} verify \eqref{s''pi}. Represent the set
$[-\pi,\pi]\cap(\mathop{\cup}\li_{j\in H_3 } I_j),$ as a union of nonintersecting
intervals $[a_\mu ,b_\mu ],\ \mu =1,...,2s+1,\ b_{\mu +1}<a_\mu .$ Let $\ul j=\ul j(\mu )$ and   $\ol
j=\ol j(\mu )$ denote the indexes $j$ such that $x_{\ul j}=a_\mu $
and $x_{\ol j}=b_\mu ,$ respectively. For each $\mu =1,...,2s+1,$
set
$$
G_\mu :=\lo(d_{\ul j+1},d_{\ol j}\rr],\quad\quad
G:=\mathop{\bigcup}\li_{\mu =1}^{2s+1}G_\mu .
$$
Without loose  of any generality verify \eqref{s''pi} only for one  interval
$G_\mu ,$ i.e., fix $\mu ,$ and let it be such that $\Pi(x)>0, \
x\in G_\mu .$ For a conveniens let $n>\ul j$ and $\ol j >3-n,$ the
cases $n=\ul j$ and  $\ol j =3-n$ are proved analogously with
respecting (d.5).

Let
$$
\oh_\mu :=\lo\{\ul j+1,...,\ol j\rr\}.
$$
Note, $\oh_\mu \subset H_3 .$ It follows from \eqref{26}, (d.0)-(d.3) that  the function $S',$ at the points $a_\nu $ defined separately for each
$\Psi_j$ with $j\in \ol H_\mu ,$ satisfies the inequality
\begin{equation}\label{28}
S'(a_\nu -)\le S'(a_\nu +).
\end{equation}

Note, $F_j\ge 0$ for $j\in\lo\{\ul j+2,...,\ol
j-1\rr\}=:\ooh_\mu\subset H_1 .$  Therefore, in particular, it follows from
the inequalities $F_{j+1}\le F_j\ge F_{j-1}$ that
\begin{equation}\label{29}
\Phi_j\,\Pi(x_j)\le 0, \quad j\in\ol H_\mu .
\end{equation}
Taking this into account, remark that in $\ol H_\mu $
there is not any $j$ for  which, in according with the definitions
(d.0)-(d.4), the following settings where made
$$
\Psi_j=\Psi_{j,3}\quad\text{and}\quad\Psi_{j-1}=\Psi_{j-1,1},
$$
as well as the settings like
$$
\Psi_{j+1}=\Psi_{j+1,3}\quad\text{and}\quad\Psi_{j}=\alpha_j\,\Psi_{j,1}+
(1-\alpha_j)\,\Psi_{j,3}\quad\text{and}\quad\Psi_{j-1}=\Psi_{j-1,1}.
$$
By the other words,
\begin{equation}\label{30}
a_\nu \ \text{(defined for}\ \Psi_j)\quad\le\quad  a_\nu\
\text{(defined for}\ \Psi_{j-1}).
\end{equation}
From this and \eqref{17} note,
\begin{equation}\label{31}
A_j''(x)=0,\quad x\notin \lo(\ul a_{j+1},\ol a_{j-1}\rr],
\end{equation}
where $\ul a_j:=a_1$ and $\ol a_j:=a_3$ if (d.3) otherwise $\ul
a_j=\ol a_j$ denote $a_\nu $ defined for $\Psi_j$ by (d.0)-(d.2) or
(d.4) (if $\ul a_{j+1}=\ol a_{j-1}$ then $\lo(\ul a_{j+1},\ol
a_{j-1}\rr]:=\emptyset $).

Using the equality  $\ul A_{j+1}=\ol A_{j},$ extract from \eqref{27}
four addends involving the function $\Psi_j$
\begin{equation}\label{32}
-F_{j+1}\,\ol A_j(x)+F_{j}\,\ol A_j(x)-F_{j}\,\ul
A_j(x)+F_{j-1}\,\ul A_j(x).
\end{equation}

Taking into account \eqref{29}-\eqref{32}, fix $j\in\ol H_\mu ,$ and show that
$$
\hskip-5.2cm\aligned (c.0) \\ (c.1) \\ (c.2) \\ (c.3)  \endaligned
\quad\quad\quad\quad\quad\quad\quad S''(x)\ge 0,\quad \lo\{
\aligned x\in(a_1,a_3]&\quad\text{if}\quad
(d.0), \\
x\in(a_1,a_2]&\quad\text{if}\quad (d.1),\\
x\in(a_2,a_3]&\quad\text{if}\quad (d.2),\\
x\in(a_1,a_3]&\quad\text{if}\quad (d.3),
\endaligned
\rr .
$$
Only these three points $a_1,\ a_2$ and $a_3$ will take part in the sentences below.

We start from the case (c.1). Describe it on $(a_1,a_2].$ The
function $\Psi_{j+1}$ can be stated by (d.0) or (d.4) or (d.1) only,
whereas $\Psi_{j-1}$ is any of the four by (d.0)-(d.4). Anyway,
$\Psi_{j+1}''=6(x-a_1),\ \ \Psi_{j}''=6(x-a_2)$ and
$\Psi_{j-1}''=0.$ Hence, by virtue of (17), write
$$
F_{j+1}\,2-F_{j+1}\,2+F_{j}\,2-F_j\frac{6(x-a_2)}{x_{j-3}-x_j}+
F_{j-1}\frac{6(x-a_2)}{x_{j-3}-x_j}\ge 0,\quad x\in (a_1,a_2],
$$
since $F_j\ge F_{j-1}.$

In the case (c.2) $\Psi_{j+1}$ is any of the four potential
settings, whereas $\Psi_{j-1}$ is  defined by (d.0) or (d.4) or (d.2)
only, but always
$$
\ol
A_j''(x)=2+\frac{6(x-a_2)}{x_{j-2}-x_{j+1}}\quad\text{and}\quad
\ul A_j''(x)=0\quad\text{for}\quad x\in(a_2,a_3],
$$
where we used (17) in the first equality. Thus,
$$
F_{j+1}\,2-F_{j+1}\,\lo(2+\frac{6(x-a_2)}{x_{j-2}-x_{j+1}}\rr)+
F_{j}\lo(2+\frac{6(x-a_2)}{x_{j-2}-x_{j+1}}\rr)\ge 0,\quad x\in (a_2,a_3],
$$
since $F_{j+1}\le F_j.$

To see (c.3) note that $\Psi_{j+1}$ and $\Psi_{j-1}$ are defined by
(d.0) or (d.4) or (d.1) and by (d.0) or (d.4) or (d.2), respectively.
Anyway, $\Psi_{j+1}''(x)=6(x-a_1),\ \
\Psi_{j}''(x)=\alpha_j\,6(x-a_2)$ and $\Psi_{j-1}''(x)=0$ for
$x\in (a_1,a_3].$ Write
$$
F_{j+1}\,2-F_{j+1}\,\frac{6(x-a_1)-\alpha_j6(x-a_2)}{3h}+
F_{j}\,\frac{6(x-a_1)-\alpha_j6(x-a_2)}{3h}
$$
$$
-F_j\,\frac{\alpha_j6(x-a_2)}{3h}+F_{j-1}\,
\frac{\alpha_j6(x-a_2)}{3h}=
F_j\lo(2+\frac{(1-\alpha_j)6(x-a_2)}{3h}-
\frac{\alpha_j6(x-a_2)}{3h}\rr)
$$
$$
+F_{j-1}\frac{\alpha_j2(x-a_2)}{h}-
F_{j+1}\frac{(1-\alpha_j)2(x-a_2)}{h}=:B_1(x)+B_2(x).
$$
So, $B_2(x)=0$ due to  the choosing of $\alpha_j$ whereas
$B_1(x)\ge 0$ for any $\alpha_j\in[0,1].$ Really, like (17)
rewrite
$$
B_1(x)=F_j\lo(2-(1-\alpha_j)\frac{6(x-a_1)-6(x-a_2)-6(x-a_1)}{3h}
\rr .
$$
$$
-\lo .
\alpha_j\frac{6(x-a_2)-6(x-a_3)+6(x-a_3)}{3h}\rr)
$$
$$
=F_j\lo(2-(1-\alpha_j)\,2+(1-\alpha_j)\frac{2(x-a_1)}{h}-
\alpha_j\,2-\alpha_j\frac{2(x-a_3)}{h}\rr)\ge 0,\quad
x\in (a_1,a_3].
$$

\medskip

For the last case (c.0) note that $\Psi_{j\pm 1}$ can both be any
of the four potential settings but it's sufficient to verify this
case only when $\Psi_{j+1}=\Psi_{j+1,2}$ and
$\Psi_{j-1}=\Psi_{j-1,2}$ because for the other settings the
positivity of $S''$ on $(a_1,a_2]\cup (a_2,a_3]$ is guaranteed
by just considered three cases, namely, on $(a_1,a_2]$ - by (c.2)
or (c.3), and on $(a_2,a_3]$ - by (c.1) or (c.3). So, for
$x\in(a_1,a_3]$ we have
$$
\ol
A_j''(x)=\frac{6(x-a_1)-6(x-a_2)_+}{3h}\quad\text{and}\quad
\ul A_j''(x)=\frac{6(x-a_2)_+}{3h},
$$
that together with (17) yield
$$
F_{j+1}\,2-F_{j+1}\,\ol A_j''(x)+ F_{j}\,\ol A_j''(x)-F_{j}\,\ul
A_j''(x)+F_{j-1}\,\ul A_j''(x)\ge 0.
$$
The inequalities (c.0)-(c.3) are proved.

Finally, since the intervals in (c.0)-(c.3) cover all $G_\mu $ if
$j$ runs through $\ol H_\mu , $  then
\begin{equation}\label{37}
S''(x)=\sum\limits_{j\in\ooh_\mu} F_j\,A_j''(x)\ge 0,\quad x\in
G_\mu ,
\end{equation}
that together with \eqref{28} leads to \eqref{s''pi}.

To prove \eqref{f-s} we need the estimate
\begin{equation}\label{38}
\lo|\Phi_j\rr|\le c \frac{\om_4(f,h)}{h^4},
\end{equation}
see, for example,
 in \cite{dzay_shev}, \eqref{1213} and the technical spline
 $$
 s(x)=L_3(x;x_n;t)+4h\sum\limits_{j=3-n}^{n-1}
\Phi_j \Psi_3(x,x_j),
 $$
 that interpolates $f$ without restrictions by cubic parabolas in each $x_j,$ see, \cite{DzGi}.
Now let $x\in [x_{j^*+1},x_{j^*-3}],$ then it follows  from \eqref{26} that
$$
\aligned
\lo|f(x)-S(x)\rr|&=\lo|f(x)- s(x)+ s(x)+S(x)\rr|\\
&\le c\,\om_4(f,h)+4h\sum\limits_{j=3-n}^{n-1}|\Phi_j|\lo|\Psi_3(x,x_j)-
\Psi_j(x)\rr|\\
&=
c\,\om_4(f,h)+\sum\limits_{j=\max\{3-n,j^*-3\}}^{\min\{n-1,j^*+3\}}
|\Phi_j|4h\lo|\Psi_3(x,x_j)- \Psi_j(x)\rr|\le
c\,\om_4(f,h)
\endaligned
$$
and therefore \eqref{f-s} is correct.

\sect{Proof of Theorem 1}

Denote the numbers
$$
b_1:=s+2,\quad b_2:=3(s+1),
$$
$$
c_9:=\max\left\{\frac{6((2\pi)^{2b_2}\max\{c_1(b_2),c_7(b_2)\}+c_8(b_2)+2)}{3c_3(b_1)},2\right\},
$$
$$
n_1:=2\left[c_9+1\right]n,\quad h_1:=h_{n_1},
$$
$$
c_{10}:=\max\left\{c_5(b_2)\left(\frac{c_8(b_2)}{2c_9}+c_1(b_2)\right),10\right\},
$$
$$
n_2:=2\left[c_{10}+1\right]n_1,\quad h_2:=h_{n_2},
$$
where $[\cdot]$ stands for the integer part.

Fix $j=3-n,...,n-1.$ For each
point $a_\nu,\ \nu =1,2,3,$ let $j_\nu$ denotes the index such that $x_{j_\nu}:=x_{j_\nu,n_1}=a_\nu,$ and let $j_\nu^*$ denotes the index such that $x_{j_\nu^*}:=x_{j_\nu^*,n_2}=x_{j_\nu} (=x_{j_\nu,n_1}).$

Let $j\in H_3.$ For each $j_\nu,\ \nu=1,2,3,$ we take
$$
{\wti{\tau}}_{j_\nu^*}(x)={\wti{\tau}}_{j_\nu^*,n_2}(x,b_2),
\quad
{\wti{t}}_{j_\nu^*}(x)={\wti{t}}_{j_\nu^*,n_2}(x,b_2),
$$
and put
\begin{align*}
\varphi_{j,\nu}(x)&:=6\int\limits_{d_j-\pi}^{x}\left({\wti{\tau}}_{j_\nu^*}(u)
+\wti h_\nu\left(\alpha_\nu{\wti{t}}_{(j_\nu+1)^*}(u)+
(1-\alpha_\nu){\wti{t}}_{(j_\nu-1)^*}(u)\right)\right)du,\quad \nu=1,3,\\
\varphi_{j,2}(x)&:=6\int\limits_{d_j-\pi}^{x}\left({\wti{\tau}}_{j_2^*}(u)
-\frac{1}{12}h^2\left(\alpha_2 t_{(j_2+5)^*}'(u)+
(1-\alpha_2)t_{(j_2-5)^*}'(u)\right)\right)du,
\end{align*}
where $\alpha_\nu\in[0,1], \nu=1,2,3,$ can be chosen such that
\begin{equation}\label{var=pi}
\varphi_{j,\nu}(d_j+\pi)=3(\pi +h)(\pi -h),\quad\nu=1,3,\quad\varphi_{j,2}(d_j+\pi)=3\pi^2-\pi h^2/2.
\end{equation}
Indeed, for example, using \eqref{x-tau}, \eqref{chi-tj} for $\wti t_{j_\nu^*},$
and \eqref{2.2'}, we, for a fixed $j,\ \nu=3$ and $\alpha_3=1,$ have the estimate
\begin{align*}
&\varphi_{j,3}(d_j+\pi)=6\int\limits_{d_j-\pi}^{d_j+\pi}\Bigl[\wti \tau_{j_3^*}(u)-(u-a_3)_+
+h\left(\wti t_{(j_3+1)^*}(u)-\chi(u,x_{j_3+1})\right)\\
&+h\bigl(\chi(u,x_{j_3+1})-\chi(u,a_3)\bigl)\Bigl]du+6\int\limits_{d_j-\pi}^{d_j+\pi}\left((u-a_3)_++h\chi(u,a_3)\right)du
\ge 3(\pi^2-h^2)+6hh_1\\
&-6\left| \int\limits_{d_j-\pi}^{d_j+\pi}\left[\wti\tau_{j_3^*}(u)-(u-a_3)_+
+h\left(\wti t_{(j_3+1)^*}(u)-\chi(u,x_{j_3+1})\right)\right]du\right|\ge 3(\pi^2-h^2)
\\
&+6hh_1
-6c_8(b_2)h_2\int\limits_{d_j-\pi}^{d_j+\pi}\Gamma_{j_3^*,n_2}^{2(b_2-s-1)}(u)du
-6c_1(b_2)h\int\limits_{d_j-\pi}^{d_j+\pi}\Gamma_{(j_3+1)^*,n_2}^{2b_2-2s-1}(u)du\\
&\ge 3(\pi^2-h^2)+6hh_1-6c_5(b_2)\bigl(c_8(b_2)h_2^2+c_1(b_2)hh_2\bigr)>3(\pi^2-h^2),
\end{align*}
whereas for $\alpha_3=0$ we (again due to $h_1>>h_2$) analogously have the opposite inequality $\varphi_{j,3}(d_j+\pi)<3(\pi^2-h^2).$ So, \eqref{var=pi} is proved for
 $\nu=3$, and for $\nu=1,2$ it can be proved by analogy.

Now, we take
$$
t_{j_\nu^*}(x)=t_{j_\nu^*,n_2}(x,b_2,Y), \quad
t_{j_\nu}(x)=t_{j_\nu,n_1}(x,b_1,Y),
$$
and put
$$
\psi_{j,\nu}(x):=\int\limits_{d_j-\pi}^{x}\Bigl[\varphi_{j,\nu}(u)
+\hat h_\nu\Bigl(\beta_\nu t_{(j_\nu+1)^*}(u)+t_{j_\nu}(u)+
(1-\beta_\nu)t_{(j_\nu-1)^*}(u)\Bigl)\Bigl]du,
$$
where $\hat h_\nu:=h^2\vee -h^2/4\vee h^2,\ \nu=1,2,3,$ respectively.

\begin{lemma}  \label{lem6}
 If a fixed $j$ belongs to $H_2,$ then $\beta_\nu\in[0,1],\ \nu=1,2,3,$ can be chosen such that
\begin{equation}\label{inpi}
\psi_{j,\nu}(d_j+\pi)=(\pi+h)\pi(\pi-h),
\end{equation}
and then three functions $\psi_{j,\nu}$ satisfy the
inequalities
\begin{equation}\label{5657}
\begin{array}{ll}
\lo(\psi_{j,\nu}''(x)-\Psi_{j,\nu}''(x)\rr)\Pi(x)\Pi(x_j)\ge
0,\\
\lo(\psi_{j,2}''(x)-\Psi_{j,2}''(x)\rr)\Pi(x)\Pi(x_j)\le
0,
\end{array}
\quad \nu=1,3,\quad x\in[-\pi,\pi] ,
\end{equation}
\begin{equation}\label{58}
\lo|\Psi_{j,\nu}(x)-\psi_{j,\nu}(x)\rr|\le c\,
h_{j,n}^3\,\Ga_{j,n}^6(x),\quad \nu=1,2,3,\quad x\in[-\pi,\pi].
\end{equation}
In addition,
\begin{equation}\label{form}
\aligned
\psi_{j,\nu}(x)=&\frac{1}{8\pi}x^4+\frac{\pi-d_j}{2\pi}x^3+
\frac{5d_j^2-6d_j\pi-h^2}{4\pi}x^2+\frac{(\pi-d_j)(5d_j^2-2\pi^2-h^2)}{2\pi}x\\
+&Q_{\jn}(x),\quad \nu=1,2,3,
\endaligned
\end{equation}
where  $Q_{\jn}\in\Bbb T_{cn_2}.$
\end{lemma}

\proof
The relations \eqref{inpi}--\eqref{58} can be proved with the arguments similar to proving \eqref{var=pi}, or \cite[Lemma 5]{Zalizko}, or \cite[Lemma 6]{dz_comono}, using the choice of $n_1$ and $n_2,$ and the inequalities
$
\Gamma_{(j_{\nu}\pm 1)^*,n_2}(x)<\Gamma_{j_{\nu}\pm 1,n_1}(x)<2\pi \Gamma_{j_{\nu},n_1}(x)<2\pi
\Gamma_{j,n}(x),\ x\in\Bbb R.
$
We will calculate here the presentation \eqref{form} only, with $\nu=1,$ for definiteness. By \eqref{tjR} and \eqref{taujR}
write
$$
\wti t_{j_1^*}(x)=\frac{1}{2\pi }x+\hat R_{j_1^*}(x),\quad
\wti\tau_{j_1^*}(x)=\frac{1}{4\pi }x^2+\frac{\pi -x_j}{2\pi}x+\tilde R_{j_1^*}(x),
$$
$$
\hat r_{j_1^*}(x):=\hat R_{j_1^*}(x)-\hat R_{j_1^*,0},\quad
\tilde r_{j_1^*}(x):=\tilde R_{j_1^*}(x)-\tilde R_{j_1^*,0},
$$
where $\hat R_{j_1^*,0}$ and $\tilde R_{j_1^*,0}$ are free terms of polynomials
$\hat R_{j_1^*}, \tilde R_{j_1^*}\in \Bbb T_{cn},$ respectively. Then
$$
\varphi_{j,1}(x)=\left(\frac{1}{2\pi}x^3+\frac{3(\pi-x_j)}{2\pi}x^2+6\tilde R_{j_1^*,0}x\right)-
\Bigl(...(d_j-\pi)\Bigl)
$$$$
-6h\left(\frac{1}{4\pi}x^2+\left(\alpha_1\hat R_{(j_1+1)^*,0}+(1-\alpha_1)\hat R_{(j_1-1)^*,0}\right)x\right)
+6h\Bigl(...(d_j-\pi)\Bigl)
$$
$$
+6\int\limits_{d_j-\pi}^{x}\Bigl(\tilde r_{j_1^*}(u)-h
\bigl(\alpha_1\hat r_{(j_1+1)^*}(u)+(1-\alpha_1)\hat r_{(j_1-1)^*}(u)\bigr)\Bigl)du
$$
$$
=\frac{1}{2\pi}x^3+\frac{3(\pi-d_j)}{2\pi}x^2+6Ax-\left(\frac{1}{2\pi}(d_j-\pi)^3+\frac{3(\pi-d_j)}{2\pi}(d_j-\pi)^2+
6A(d_j-\pi)\right)
+q_{j_1}(x),
$$
where
$$
A:=\tilde R_{j_1^*,0}-h\left(\alpha\hat R_{(j_1+1)^*,0}+(1-\alpha)\hat R_{(j_1-1)^*,0}\right),
$$
and $q_{j_1}\in\Bbb T_{cn}$ does not have a free term.
Taking this and \eqref{var=pi} we derive the value of $A$ 
$$
3(\pi^2-h^2)=\frac{1}{2\pi}\bigl((d_j+\pi)^3-(d_j-\pi)^3\bigr)+
\frac{3(\pi-d_j)}{2\pi}\bigl((d_j+\pi)^2-(d_j-\pi)^2\bigr)+12\pi A
$$
$$
\Rightarrow \ A=\frac{5d_j^2-6d_j\pi-3h^2}{12\pi},
$$
and so,
$$
\aligned
\varphi_{j,1}(x)=&\frac{1}{2\pi}x^3+\frac{3(\pi-d_j)}{2\pi}x^2+\frac{5d_j^2-6d_j\pi-3h^2}{2\pi}x+ \frac{(\pi-d_j)
(3d_j^2-2d_j\pi-2\pi^2-3h^2)}{2\pi}\\
+&q_{j_1}(x).
\endaligned
$$
Having this, \eqref{tjR} and \eqref{inpi} we get \eqref{form} analogously. Lemma \ref{lem6} is proved.

\eop

\vskip 0.3cm

\centerline{\it Construction of the nearly coconvex polynomial}

\medskip

 For each $j=3-n,n-1$ introduce the polynomial
$\psi_j(x)\in\Bbb T_{cn_1}.$ If $j\in H_2$ then set
$$
\psi_j(x):=\psi_{j,2}(x)\quad\text{if}\quad \Phi_j\,\Pi(x_j)\le
0,
$$
{\it otherwise} set
$$
\psi_j(x):=\lo\{ \aligned \psi_{j,1}(x)&\quad\text{if}\quad
|F_{j+1}|>|F_{j}|\ge |F_{j-1}|,\\
\psi_{j,3}(x)&\quad\text{if}\quad |F_{j+1}|\le |F_{j}|<
|F_{j-1}|,\\
\alpha_j\,\psi_{j,1}(x)+(1-\alpha_j)\,\psi_{j,3}(x)&\quad\text{if}\quad
|F_{j+1}|>|F_{j}|< |F_{j-1}|.
\endaligned
\rr .
$$
If $j\notin H_2$ (i.e., $j:\ x_j\in O_{i,2},\ i=1,...,2s,$)
then let
$$
\psi_j(x):=\lo\{ \aligned
\psi_{j,2}(x)&\quad\text{if}\quad\Phi_j\,\Pi(x_j,\ti Y_i)\le 0,\\
\psi_{j,1}(x)&\quad\text{otherwise.}
\endaligned \rr .
$$

Now, put
\begin{equation}\label{p}
P_n(x)=L_3(x;x_n;f)+4h\sum\limits_{j=3-n}^{n-1}\Phi_j\,
\psi_j(x).
\end{equation}\label{p}
The fact that $P_n$ is a polynomial from $\Bbb T_{cn}$ can be directly verified arithmetically like in 
\cite{Zalizko}, or \cite{dz_comono}, using \eqref{form}, i.e., all the arithmetical terms in \eqref{form}, 
having been evaluated in the sum \eqref{p} together with the corresponding divided differences, including the $L_3,$ 
are equal $0.$

Verify \eqref{nearly}. Remark that Lemma \ref{lem6} will be used in two senses: in an
"ordinary" one for $j\in H_2 =H (n,Y,2),$ and for $j\notin H_2$ in the sense that $j\in H (n,\ti Y_i,2).$
So, \eqref{5657}, \eqref{14},
\eqref{26}, \eqref{27} and \eqref{37} imply
$$
P_n''(x)\Pi(x)=\lo(
L_3''(x;x_n;f)+4h\sum\limits_{j=3-n}^{n-1}\Phi_j
\bl \psi_j''(x)- \Psi_j''(x)\br +4h\sum\limits_{j=3-n}^{n-1}\Phi_j
\Psi_j''(x)\rr) \Pi(x)
$$
$$
=4h\sum\limits_{j\in H_2}\frac{1}{\Pi^2(x_j)}\Phi_j\Pi(x_j)
\lo( \psi_{j,\nu}''(x)- \Psi_{j,\nu}''(x)\rr)\Pi(x)\Pi(x_j)
$$
$$
+4h\sum\limits_{i=1}^{2s}\sum\limits_{j:x_j\in
O_{i,1}}\frac{1}{\Pi^2(x_j,\ti Y_i)}\Phi_j\Pi(x_j,\ti Y_i)
\lo( \psi_{j,2\vee 1}''(x)- \Psi_{j,2\vee
1}''(x)\rr)\Pi(x,Y)\Pi(x_j,\ti Y_i)
$$
$$
+\lo(F_n\bbl
2-\frac{\Psi_n''(x)-\Psi_{n-1}''(x)}{3h}\bbr
+\sum\limits_{j=3-n}^{n-1}F_j\,A_j''(x)+F_2\frac{\Psi_3''(x)}{3h}\rr)\Pi(x)
=:A(x)+B(x)+C(x),
$$
$$\aligned
A(x)\ge 0,&\quad\quad\quad x\in\Bbb R,\\
B(x)\ge 0,&\quad\quad\quad x\in\Bbb R\setminus\cup_{i\in\Bbb Z}(x_{j_i+5},y_i),\\
C(x)\ge 0,&\quad\quad\quad x\in G \ \text{on all periods},
\endaligned
$$
that leads to \eqref{nearly}. 
To prove \eqref{estimate} we use \eqref{f-s}, \eqref{38}, \eqref{58} and
\eqref{2.2}. Namely,
$$
\lo\|f-P_n\rr\|=\lo\|f-S+S-P_n\rr\|
=\lo\|f-S+\sum_{j=3-n}^{n-1} \Phi_j \, 4h\,\bl
\Psi_j(\cdot) -\psi_j(\cdot)\br\rr\|_{[-\pi,\pi]}
$$$$
\le c\,\om_4(f,h)+c\,\lo\|\sum_{j=3-n}^{n-1}\om_4(f,h)\Ga_j^6(\cdot)\rr\|_{[-\pi,\pi]}\le c\,\om_4(f,h).
$$
Theorem 1 is proved. 
\eop



\renewcommand\baselinestretch{0.5}

\end{document}